\documentclass[psamsfonts]{amsart}

\usepackage{latexsym}
\usepackage{amscd}
\usepackage{amssymb}
\usepackage{amsmath}
\usepackage{xy}
\xyoption{all}
\usepackage{amsthm}
\usepackage{enumerate}
\usepackage{hyperref}
\usepackage{pstricks}
\usepackage{epsfig}
\usepackage{graphics}

\newtheorem{thm}{Theorem}[section]
\newtheorem{lem}[thm]{Lemma}
\newtheorem{defn}[thm]{Definition}
\newtheorem{prop}[thm]{Proposition}
\newtheorem{cor}[thm]{Corollary}
\newtheorem{eg}[thm]{Example}

\newtheorem{rmk}[thm]{Remark}

\numberwithin{equation}{section}

\def\C{\mathbb{C}}

\def\Z{\mathbb{Z}}

\def\P{\mathbb{P}}
\def\O{\mathcal{O}}
\def\del{\partial}
\def\delbar{\overline{\partial}}

\def\ra{\rightarrow}

\def\del{\partial}
\def\delbar{\overline{\partial}}
\def\Ainf{A_{\infty}}

\def\Fuk{\operatorname{Fuk}}
\def\Qcoh{\operatorname{Qcoh}}
\def\Coh{\operatorname{Coh}}

\def\Func{\operatorname{Func}}
\def\im{\operatorname{im}}
\def\mod{\operatorname{mod}}
\def\Mod{\operatorname{Mod}}
\def\dgMod{\operatorname{dgMod}}

\def\perf{\operatorname{perf}}

\def\id{\operatorname{id}}

\def\Hom{\operatorname{Hom}}

\def\Ext{\operatorname{Ext}}
\def\Sym{\operatorname{Sym}}

\def\rk{\operatorname{rank}}
\def\tr{\operatorname{tr}}

\def\ker{\operatorname{ker}}
\def\coker{\operatorname{cok}}
\def\charac{\operatorname{char}}
\def\Spec{\operatorname{Spec}}

\begin{document}

\author{Matthew Robert Ballard}
\address{Department of Mathematics, University of Washington,
Seattle, WA 98195, USA}
\email{ballard@math.washington.edu}
\title{Sheaves on Local Calabi-Yau Varieties}
\maketitle

\begin{abstract}
 We investigate sheaves supported on the zero section of the total space of a locally-free sheaf $E$ on a smooth, projective variety $X$ when $E$ satisfies $\bigwedge^{\rk E} E \cong \omega_X$. We rephrase this construction using the language of $\Ainf$-algebra and provide a simple characterisation of the case $E = \omega_X$.
\end{abstract}

\section{Introduction}

Calabi-Yau varieties formed by taking the total space of a locally-free sheaf $E$ on a variety $X$ with $\bigwedge^d E \cong \omega_X$ are important testing ground for ideas in mathematics and physics. They are often viewed as an approximation to a proper Calabi-Yau variety. In this paper, we study these local Calabi-Yau varieties from a homological perspective. We are interested in the bounded derived category of coherent sheaves with support on the zero section of $V(E)$, $D^b_X(\Coh(V(E)))$, and the larger unbounded derived category of quasi-coherent sheaves $D_X(\Qcoh(V(E)))$. In the case that $X$ is smooth and proper, $D^b_X(\Coh(V(E)))$ possesses a trivial Serre functor, hence provide simple examples of Calabi-Yau triangulated categories. Then, we have two natural maps $i: X \hookrightarrow V(E)$ and $\pi: V(E) \ra X$ satisfying $\pi \circ i = \id_X$. These furnish functors which descend to the derived categories $i_*,\pi^*: D(\Qcoh(X)) \ra D(\Qcoh(V(E)))$ and $Li^*,\pi_*: D(\Qcoh(V(E))) \ra D(\Qcoh(X))$. The image of $i_*$ lies in $D_X(\Qcoh(V(E)))$ and generates it as a triangulated category. The functor $i_*$ is not full but as is well-known \cite{ST00}
\begin{displaymath}
 \Ext^j_E(i_*F,i_*G) = \bigoplus_{r=0}^{\rk E} \Ext^{j+r}_X(F,G \otimes \bigwedge^r E)
\end{displaymath}
After a brief discussion of the choices of dg-models underlying $D^b(\Coh(X))$ and of sheaves with support on sub-schemes, we lift this isomorphism to a quasi-isomorphism on the chain level and give the proper algebra structure on the dg-endomorphisms of $i_*F$ for this quasi-isomorphism to be a morphism of dg-algebras. This allows us to describe $D_X(\Qcoh(V(E)))$ and $D_X^b(\Coh(V(E)))$ in simple terms of $X$ and $E$. We can cut the data necessary by appealing to the triviality of Serre duality on $D^b_X(\Coh(V(E)))$. At the chain level, this manifests itself as symmetric pairing on the Hom-spaces of the dg-category which is non-degenerate on cohomology. In the special case where $E = \omega_X$, one expects that we can describe $D_X(\Qcoh(\omega_X))$ solely in terms $D(\Qcoh(X))$. Indeed, this is case. The idea, while straightforward, requires some technical development to make manifest. This leads into $\Ainf$-algebras and formal symplectic geometry, which occupies the final two sections of the paper. We conclude the paper by making connections with homological mirror symmetry.

\section{DG-models for the bounded derived category of coherent sheaves}

Our main interest in this paper is not in the triangulated category
$D^b(\Coh(X))$ but in an enhanced version of it \cite{BK91}. For any
separated scheme $X$ there is a standard enhancement - namely the dg-category
of bounded below complexes of injective sheaves with bounded cohomology.
The utility is its universality. Its main drawback is its inaccessibility to
computation. We have in mind another dg-category $\check{C}(X)$.

\begin{defn}
 Let $\check{C}(X)$ denote the dg-category whose objects are locally-free
coherent sheaves and whose morphism spaces are given by choosing some
affine cover $\mathcal{U}$ and setting
\begin{equation*}
 \Hom_{\check{C}(X)}(E,F) = \check{C}(\mathcal{U},E^{\vee} \otimes F)
\end{equation*}
\end{defn}

If $X$ is a smooth projective variety over $\C$, there is another more
differential geometric construction. Although it will not be featured heavily
in this paper, it has one nice quality.

\begin{defn}
 Let $\delbar(X)$ denote the dg-category with objects locally-free coherent sheaves
and morphisms give by
\begin{equation*}
 \Hom_{\delbar(X)}(E,F) = \Gamma(V(F)^{\vee}\otimes V(E) \otimes \Omega^{0,*})
\end{equation*}
The differential is the Dolbeault differential and $V(E)$ is the total space of the geometric vector bundle corresponding to the locally-free sheaf $E$. $\Omega^{0,*}$ the bundle of $d\bar{z}$-forms.
\end{defn}

\begin{rmk}
 The main attractive quality of $\delbar(X)$ is the presentation of Serre
duality. The trace morphism $\tr: H^n(\omega_X) \ra \C$ is naturally defined. $\Gamma(\omega_X \otimes \Omega^{0,n})$
equals $\Gamma(\Omega^{n,n})$ and, hence, any section gives a top-degree differential
form $X$ which can be integrated. Applying the standard argument using
bump functions, one sees that the pairing is non-degenerate
at the chain level.
\end{rmk}

Given two locally-free sheaves $E$ and $F$, in each dg-model, the cohomology
of each morphism chain complex computes $\Ext_X^*(E,F)$. In nice situations, one therefore expects
that the derived categories are equivalent as triangulated categories.

Let us recall how to derive a dg-category. Let $\mathcal{C}$ be a
dg-category over $k$. Then, we can consider the category of dg-functors $M:
\mathcal{C}^{op} \ra Ch(k)$ where $Ch(k)$ is the dg-category of chain complexes
over $k$. Such a functor is called a right $\mathcal{C}$-module. The dg-category
$\Func(\mathcal{C}^{op},Ch(k))$ is called the category of right $\mathcal{C}$-
modules and is denoted by $\Mod \mathcal{C}$. $\Mod \mathcal{C}$ naturally has
cones and a shift functor that provide it with a pre-triangulated structure. Therefore,
$H^0(\Mod \mathcal{C})$ is a triangulated category with triangles
induced by these cones. We consider the smallest ore-triangulated and idmepotent-closed triangulated category containing the image of $\mathcal{C}$ and closed under idempotent splittings. Localising this at quasi-isomorphisms yields the (perfect) derived category of $\mathcal{C}$ denoted by $D(\mathcal{C})$. (This is usually denoted $D_{\perf}(C)$). For more information and references see \cite{Kel06}.

Let $I(X)$ denote the dg-category consisting of a choice of K-injective resolution of each (unbounded) complex of quasi-coherent sheaves, see \cite{Spa88}. Then, if $X$ is quasi-compact, $H^0(I(X)) \cong D(\Qcoh(X))$. Let us recall the following result from \cite{Sei03}.

\begin{prop}
 $D(\check{C}(X))$ is equivalent to the smallest triangulated category of $D(I(X))$ containing all locally-free coherent sheaves.
\end{prop}

\proof Given such a locally-free sheaves $E$, let $I_E$ denote a choice of injective resolution. We shall form a new category $T(X)$ whose objects are locally-free sheaves and whose morphisms are given by
\begin{equation*}
 \Hom_{T(X)}(F,E) = \begin{pmatrix} \Hom_{\mathcal{I}(X)}(I_F,I_E) & 
\check{C}(F^* \otimes I_E)[1] \\ 0 & \check{C}(F^* \otimes E) \end{pmatrix}
\end{equation*}
Consider the double complex formed by $\check{C}(F^* \otimes I_E)$. From
this we have natural homomorphisms
\begin{equation*}
\Hom_{I(X)}(I_F,I_E) \ra \check{C}(F^* \otimes I_E) \leftarrow
\check{C}(F^* \otimes E)
\end{equation*}
The first map comes from composing with the map $F \ra I_F$ and then restricting
to the open affines. The second map comes from composing with the map $E \ra I_E$.
Using these maps we can appropriately define composition in $T(X)$
\begin{equation*}
\Hom_{\mathcal{I}(X)}(I_F,I_G) \otimes \check{C}(E^* \otimes I_F) \ra \check{C}
(E^* \otimes I_G)
\end{equation*}
\begin{equation*} 
\check{C}(F^* \otimes I_G) \otimes \check{C}(E^* \otimes F) \ra \check{C}(E^* 
\otimes I_G) 
\end{equation*}
Consequently using the natural two natural maps above, we obtain functors
$I(X) \leftarrow T(X) \ra \check{C}(X)$. We wish to show that each of these is
a quasi-isomorphism. The complex $\check{C}(\mathcal{U},F^{\vee} \otimes I_E)$ is
naturally bi-graded with the decomposition of the differential into a \v{C}ech
piece and $d_{F^{\vee}} \otimes 1 + 1 \otimes d_{I_E}$. If we take cohomology
with respect the \v{C}ech differential, we will get only get the complex
$\Gamma(F^{\vee} \otimes I_E)$ since $F^{\vee} \otimes I_E$ is still a complex
of injectives. Thus, the cohomology of $\check{C}(\mathcal{U},F^{\vee} \otimes I_E)$
is simply $\Ext(F,E)$. Each map preserves the natural bi-gradings involved.
In the case of $\Hom(I_F,I_E)$ the bi-grading is trivial (i.e. is a regular grading).
The computation above then shows that the first map is a quasi-isomorphism. For
the second, we note that $F^{\vee}|_U$ for any affine $U$ is a bounded complex of
projectives. Hence, $H(\Hom(F|_U,(I_E)_U)) = H(\Hom(F|_U,E_U)$ as $E \ra I_E$ is
a quasi-isomorphism. Thus, taking the differential with respect to the other piece
produces an isomorphism on the $E^1$-page and hence induces a quasi-isomorphism.
Thus, the natural functors above induce a quasi-equivalence between $\check{C}(X)$ and the sub-category $I(X)$ consisting of resolutions of locally-free sheaves. The final statement is now clear. \qed

\begin{cor}
 If $X$ is smooth and projective over $k$, then $D(\check{C}(X))$ is triangle equivalent to $D^b(\Coh(X))$.
\end{cor}

Next we essentially repeat the argument to show that $\check{C}(X)$ and $\delbar(X)$
are quasi-equivalent.

\begin{prop}
 $D(\check{C}(X))$ and $D(\delbar(X))$ are equivalent.
\end{prop}

\proof We will combine \v{C}ech and Dolbeault cohomology into a double
complex $C^{s,t} = \check{C}^s(V(F) \otimes V(E)^{\vee} \otimes \Omega^{0,t})$ with the sum of the two
differentials. Now we proceed as in the previous proof by forming a new category $\mathcal{T}'$ whose objects are locally-free sheaves and whose morphism space between $E$ and $F$ is $\check{C}^s(V(F) \otimes V(E)^{\vee} \otimes \Omega^{0,t})$. We again have natural dg-functors from $\check{C}(X)$ and $\delbar(X)$ to $\mathcal{T}'$. These are quasi-isomorphisms as again can be checked by looking at the associated spectral sequence to $C^{s,t}$. \qed

We shall use $I(X)$ for most of the formal results in the next section and $\check{C}(X)$ when considering local Calabi-Yau varieties.

\section{Preliminaries on coherent sheaves supported on subvarieties}

Our objects of interest are coherent sheaves on $V(E)$ supported on the zero
section. In this section, we recall the general notions of categories of sheaves with fixed support. Let $Y$ be a
projective scheme over $k$ and $Z$ a closed subscheme with ideal sheaf $I_Z$.

\begin{defn}
A coherent sheaf $F$ on $Y$ is supported on $Z$ if some power of $I_Z$ annihilates $F$. A quasi-coherent sheaf is supported on $Z$ if all coherent subsheaves are supported on $Z$.
\end{defn}

The following is classical.

\begin{lem}
 If $F$ is a quasi-coherent sheaf supported $Z$, then there exists an injective quasi-coherent sheaf $I(F)$ supported on $Z$ and an injection $F \hookrightarrow I(F)$.
\end{lem}

\proof First, note that for any point $x \not \in Z$, we have $F_x = 0$. To get the desired
quasi-coherent sheaf, choose the injective envelopes $I(F_x)$ of $F_x$ for all $x \in Z$
and let $I(F) = \prod_{x\in Z} i_x^*(I(F_x))$ where $i_x$ is the inclusion of the point in
$Y$. Now, let $J = \{i \in I(F_x) | (I_Z)_x^ki = 0 \text{ for some } k\}$. Then, $J$ is
injective and $F_x \subset J$. Hence, $J = I(F_x)$. $I(F)$ is therefore
supported on $Z$ and clearly $F$ injects into $I(F)$. \qed

Let $D_Z(\Qcoh(Y))$ denote the subcategory of $D(\Qcoh(Y))$ consisting of complexes whose cohomology sheaves are supported on the zero section.

\begin{lem}
 $D(\Qcoh_Z(Y))$ is triangle equivalent to $D_Z(\Qcoh(Y))$.
\end{lem}

\proof %Clearly there is an inclusion of $D(\Qcoh_Z(Y))$ into $D_Z(\Qcoh(Y))$. The previous lemma tells us that there are enough injectives in $\Qcoh_Z(Y)$ which remain injective in $\Qcoh(Y)$. Thus, if we take the K-injective resolution on a complex of sheaves supported on the $Z$, it is K-injective in $Ch(\Qcoh(Y))$. Thus, the inclusion is indeed full and faithful. We just need to check its image. Now let us look at $\leftexp{\perp}{D(\Qcoh_Z(Y))}$ which is defined as the category of all objects $C$ of $D(\Qcoh(Y))$ for which $\Hom_{D(\Qcoh(Y))}(C,D)$ is zero for all $D$ in $D(\Qcoh_Z(Y))$. No nonzero object of $D_Z(\Qcoh(Y))$ lies in $\leftexp{\perp}{D(\Qcoh_Z(Y))}$. For an complex $C$ with a nonzero cohomology sheaf supported on $Z$ in degree $j$, take the kernel $\ker d_j$ and consider the map constructed in the previous lemma by localising at each point of $x \in Z$. Since the cohomology sheaf is supported on $Z$ and nonzero, this is the desired nonzero morphism. Now, for a given $C$ from $D_Z(\Qcoh(Y))$, one constructs the homotopy limit $\hocolim_Z(C)$ over maps from $C$ to $D(\Qcoh_Z(Y))$. (More precisely, we take to be compact). This sits in a triangle
% \begin{center}
% \leavevmode
% \begin{xy}
%  (0,5)*+{C}="a"; (10,-10)*+{\hocolim_Z(C)}="b"; (-10,-10)*+{C(i)}="c"; {\ar@{->}^i "a";"b"}; {\ar@{->} "b";"c"}; {\ar@{->}^{[1]} "c";"a"}
% \end{xy}
% \end{center}
% $C(i)$ is still in $D_Z(\Qcoh(X))$ since it is a triangulated subcategory. But, from the construction of the homotopy colimit, $C(i) \in \leftexp{\perp}{D(\Qcoh_Z(Y))}$. Thus, $C(i) \cong 0$ and $C$ is isomorphic to an object of $D(\Qcoh_Z(Y))$.

It is enough to prove this result for bounded below complexes of injectives in $D_Z(\Qcoh(X))$ since we use these to build K-injective resolutions. Let $I_*$ be a bounded below complex of injective sheaves whose cohomology sheaves are supported on $Z$. We shall construct a complex of injective sheaves $J_*$ supported on $Z$ quasi-isomorphic to $I_*$. Assume that the first non-zero term in $I$ is $I_0$. Then, $\ker d_0 \subset I_0$ is a coherent sheaf supported on $Z$. From the previous lemma, we can choose an injective sheaf $J_0$ and an injection $\ker d_0 \hookrightarrow J_0$. This extends to a morphism $\psi_0:I_0 \ra J_0$. $\coker \psi_0$ is also supported on $Z$. Choose an injective
sheaf $J_1'$ supported on $Z$ and so that $\coker \psi_0$ injects into $J_1'$. The
induced map $I_0 \ra J_1'$ factors through $\im d_0$ by construction, hence induces a map $\psi_1' : I_1 \ra J_1'$. Choose an injective sheaf $J_1''$ supported on $Z$
so that $\ker d_1 / \coker d_0$ injects into $J_1''$. Set $J_1 = J_1' \oplus J_1''$ with the differential $d_0: J_0 \ra J_1$ given by the map to $J_1'$. Then we have an induced map $\psi_1: I_1 \ra J_1$, $\psi_0$ is an isomorphism on $H^0$, $\psi_1$ is an injection on $H^1$. Now iterate this procedure to give $J_*$. \qed

\begin{cor}
  $D^b(\Coh_Z(Y))$ is equivalent to $D^b_Z(\Coh(Y))$.
\end{cor}

\proof Assume that $J_*$ is a bounded below complex of injectives with support on $Z$ and that $J_*$ has bounded coherent cohomology. We need to show that $J_*$ is quasi-isomorphic to a bounded complex of coherent sheaves with support on $Z$. To do this we let $K_*$ and $B_*$ be the kernel and images of the differentials. We proceed again by induction and assume that $J_*$ begins at zero. We take $E_0 = K_0$. To construct $E_1$ we note that since $H^1(J_*)$ is coherent we can find a coherent subsheaf $F_1 \subset K_1$ which surjects onto $H^1(J_*)$. The kernel of the map $F_1 \ra H^1(J_*)$ lies in the image of $d_0$ and is coherent. Thus, we can redefine $E_0$ as the preimage of this kernel under $d_0$ which is coherent. Now we proceed by induction, noting that to finish the $i$-the stage we only need to redefine $E_{i-1}$. \qed

\begin{rmk}
 This answers part of a seemingly unanswered question from \cite{Bri05}, see the paragraph preceding lemma $4.4$.
\end{rmk}

\begin{lem}
 Let $i: Z \ra Y$ denote the inclusion and $i_*: \Coh(Z) \ra \Coh(Y)$ ($i_*: \Qcoh(Z) \ra \Qcoh(Y)$) the induced functor. Then, $i_*$ is exact and the image of $i_*$ generates $\Coh_Z(Y)$ as an abelian category, i.e. the smallest abelian subcategory of $\Coh_Z(Y)$ containing the essential image of $i_*$ is $\Coh_Z(Y)$ itself. Similarly, the smallest abelian subcategory closed under arbitrary direct sums containing the image of $i_*$ in $\Qcoh_Z(Y)$ is $\Qcoh_Z(Y)$ itself. 
\end{lem}

\proof $i_*$ is always left-exact and since it is an isomorphism onto its image
it is right exact. Any object $F$ in $\Coh_Z(Y)$ admits a filtration $0 \subset F_0
\subset F_1 \subset \cdots \subset F_n \subset F$ such that $F_i/F_{i-1}$ is
annihilated by $I_Z$. Thus, it is of the form $i_*E$ for some $E$. The resulting
short-exact sequences show that we can get any object of $\Coh_Z(Y)$ using a finite number (which depends on $F$ and is unbounded over $\Coh_Z(Y)$) of iterations of short exact sequences starting from objects in the image of $i_*$. The final statement results from the fact any quasi-coherent sheaf is a union of its coherent subsheaves. \qed

We say that a subcategory $\mathcal{S}$ of a triangulated category $\mathcal{T}$ strongly generates if the smallest triangulated category of $\mathcal{T}$ containing $\mathcal{S}$ is $\mathcal{T}$ itself. We saw that it generates if the smallest triangulated category of $\mathcal{T}$ closed under direct sums and containing $\mathcal{S}$ is $\mathcal{T}$ itself. The following in now immediate from the previous lemma.

\begin{cor}
 The image of $D^b(\Coh(Z))$ under $i_*$ strongly generates $D^b(\Coh_Z(Y))$ and the image of $D(\Qcoh(Z))$ under $i_*$ generates $D(\Qcoh_Z(Y)$.
 \label{cor:generationcor}
\end{cor}

\begin{rmk}
 In \cite{Rou03}, a notion of dimension of a triangulated category is defined. It is noted that the dimension of $D^b_Z(\Coh(Y))$ is infinite. Indeed, any complex must be annihilated by some large power of $I_Z$. Thus, an infinite number of extensions is required reach all sheaves supported on $Z$. Consequently, $D^b_Z(\Coh(Y))$ cannot be smooth in the sense of \cite{KS06} although it is compact.
\end{rmk}

Often one wants to view these local varieties as small neighborhoods of $Z$ in some
larger ambient variety. In the Zariski topology, this of course is problematic. To
make it precise we must pass to the formal completion of $Y$ along $Z$.

Let us first recall the local situation. Let $R$ be a Noetherian ring and $I$ an ideal. Then the formal completion $\hat{R}$ along $I$ is the inverse limit of the system $R/I^n \ra R/I^{n-1}$. Similarly, given a module $M$ over $R$, we can complete $M$ using the inverse system $M/I^n \ra M/I^{n-1}$ to get a module $\hat{M}$ over $\hat{R}$. The following lemma is standard
\begin{lem}
 If $M$ is finitely-generated, then $\hat{M} \cong M \otimes_R \hat{R}$.
\end{lem}

Analogous to the previous discussion we can define the abelian category of finitely-
generated $\hat{R}$-modules supported along $\hat{I}$, $\mod_{\hat{I}}(\hat{R})$ or
the abelian category of modules supported along $\hat{I}$, $\Mod_{\hat{I}}(\hat{R})$. (Recall that the completion of a Noetherian ring is
again Noetherian).

\begin{lem}
 $\Mod_I(R)$ is equivalent to $\Mod_{\hat{I}}(\hat{R})$ as an abelian category.
\end{lem}

\proof Recall that $\hat{R}/\hat{I}^n \cong R/I^n$. Now consider the completion
functor from $\Mod(R)$ to $\Mod(\hat{R})$. Restricted to $\mod(R)$, it is exact.
In addition, for $M \in \mod_I(R)$ we have $\hat{M} \cong M \otimes_R \hat{R} \cong
M \otimes_{R/I^n(M)} \hat{R}/\hat{I}^n(M) \cong M$. Thus, $\mod_I(R)$ is equivalent
to $\mod_{\hat{I}}(\hat{R})$. Any element of $\Mod_I(R)$ naturally has the structure of an $\hat{R}$-module. An $\hat{r} = \lbrace r_l \rbrace \in R/I^l$ acts on $n \in N$ by $\hat{r} n = r_l n$ for $l$ so that $I^l n = 0$. Since $I^j n = 0$ for $j \geq l$, $r_j$ acts the same as $r_l$. Let $N$ be an element of $\Mod_{\hat{I}}(\hat{R})$. $N = \bigcup_{j \in J} F_j$ for $F_j$ finitely-generated for $\hat{R}$-modules. Thus, $F_j$ is annihilated by $\hat{I}^k$ for some $k$ and has the structure of an $R$-module. Consequently, $N$ is an $R$-module. \qed

This is simply the local case of completing a Noetherian scheme $X$ along a closed
subscheme $Y$. Consequently, we have
the following corollaries.

\begin{cor}
 $\Qcoh_Y(X) \cong \Qcoh_{\hat{Y}}(\hat{X})$
\end{cor}

\begin{cor}
 $D_Y(\Qcoh(X))$ and $D^b_Y(\Coh(X))$ only depend on the isomorphism class of $\hat{X}$.
\end{cor}

Let us now restrict our attention to the case where $X$ and $Y$ are smooth. As before let $I_Y$ denote the ideal sheaf of $Y$ on $X$. Then $I_Y/I^2_Y$ is the conormal sheaf of $Y$ in $X$. Recall the following fact.

\begin{lem}
 With the hypotheses above, $I_Y^n/I_Y^{n+1} \cong \Sym^n(I_Y/I^2_Y)$.
\end{lem}

We now have extensions.
\begin{equation*}
0 \ra I_Y^n/I_Y^{n+1} \ra \mathcal{O}_X/I_Y^{n+1} \ra \mathcal{O}_X/I_Y^n \ra 0
\end{equation*}
This gives us a way to determine the isomorphism class from the infinite sequence of extensions lying in $\Ext_Y^1(\mathcal{O}_X/I_Y^n, I_Y^n/I_Y^{n+1})$. There are some situations where we only need to know the sub-variety and its conormal bundle.

\begin{lem}
 If $Y$ satisfies $H^1(I_Y^n/I_Y^{n+1}\otimes(I_Y^r/I_Y^{r+1})^{\vee}) = 0$ for any
$n > r \geq 0$, then $\hat{X}$ only depends on the isomorphism class of $Y$ and it's conormal bundle in $X$.
\label{cor:completionlemma}
\end{lem}

\proof We simply compute by induction. As mentioned, $\mathcal{O}_X/I_Y^2$ is an extension of $\mathcal{O}_Y$ by $I_Y/I_Y^2$ which must be trivial if $H^1(I_Y/I_Y^2) = 0$. Iterating we see that the given condition on the cohomology guarantees that all the extensions will be trivial \qed

In particular, we have the following useful case thanks to Kodaira vanishing.

\begin{cor}
Assume the characteristic of our field is zero. If the conormal bundle of $Y$ in $X$ is the dual of the canonical bundle and $Y$ is Fano, then the formal neighborhood of $Y$ and hence $D^b_Y(\Coh(X))$ is uniquely determined.
\end{cor}

In general, we shall say that the formal neighborhood is trivial if all the extensions mentioned above vanish. The trivial formal neighborhood is isomorphic to the formal neighborhood of $Y$ in the total space of its normal bundle.

\begin{cor}
 Assuming the formal neighborhood is trivial, the Grothendieck groups $K(Y)$ and $K_Y(X) = K(D^b_Y(\Coh(X)))$ are isomorphic.
\end{cor}

\proof We have already seen that $i_*(D^b(\Coh(Y))$ generates $D^b_Y(\Coh(X))$.
We can assume we are working on total space the normal bundle $N$. Here we have $R\pi_* \circ i_* = \pi_* \circ i_* = \id$. Thus, $i_*$ is injective. \qed

\begin{rmk}
 While the Grothendieck groups may coincide, the Euler pairings definitely do not.
\end{rmk}

\section{Sheaves on local varieties}

Recall the construction of a geometric vector bundle $V(E)$ from a locally-free coherent sheaf $E$ on a scheme $X$. Choose an affine covering $U_i=\Spec A_i$ for which $E$ is free. Then $E|_{U_i} \cong A_i^n$ with a basis we shall denote by $x_i$. Consider the algebra $B_i = \Sym^*(E|_{U_i})$. Then, $V(E)$ is the scheme formed by gluing together the copies of $\Spec B_i$ using the restriction morphism $U_i,U_j \ra U_i \cap U_j$. 

One common use of this construction is to make Calabi-Yau varieties. Assume that $X$ is smooth for the next lemma.

\begin{lem}
$V(\omega_X)$ has trivial canonical bundle.
\end{lem}

\proof We will first provide a local description for a non-vanishing section of the canonical
bundle of $V(\omega_X)$. Let $U$ be an affine chart on $X$ and take a non-vanishing section $\sigma$ of $\omega_X$ over $U$. Then, $V(\omega_X)(\pi^{-1}U) = \Spec \O(U)[\sigma]$. Identifying the
sheaf of sections with $\omega_X^{-1}$ we get an element $\sigma \wedge \sigma^{-1}$ of
$\omega_{V(\omega_X)}$ over $\pi^{-1}(U)$. Any other choice of section $\sigma'$ will differ
from $\sigma$ by a non-vanishing function that will be canceled out after wedging with
$\sigma^{-1}$. Thus, this construction is independent of the choice of $\sigma$ and gives
a trivialisation of $\omega_{V(\omega_X)}$. \qed

The proof of the lemma extends to following case.

\begin{lem} Let $E$ be a locally-free coherent sheaf on $X$ of rank $r$ such that
$\bigwedge^r E \cong \omega_X$. Denote by $V(E)$ the total space of the vector bundle
associated to $E$. Then, $V(E)$ is Calabi-Yau.
\end{lem}

\begin{rmk}
 \begin{enumerate}
  \item We can require less than smoothness of $X$, \cite{Kaw02}.
  \item There is an obvious strong relation between $V(\omega_X)$ and divisors in the anti-canonical class
of $X$. Namely, sections of $V(\omega_X)$ are in bijection with said divisors. Similarly, if we
took $E = L_1 \oplus \cdots L_n$ such that $L_1 \otimes \cdots \otimes L_n = \omega_X$,
sections of $V(E)$ are in bijection with complete intersections determined by divisors
in the classes of the $L^{-1}_i$.
 \end{enumerate}
\end{rmk}

In this section, we work with the \v{C}ech dg-model for the bounded derived category of coherent sheaves on a variety.
Consider the pullback of $E^{\vee}$ to $V(E)$. The maps $x^{\vee}_i \otimes x_j \ra \delta_{ij}$ give maps $\pi^*E^{\vee}|_{\pi^{-1}(U_i)} \ra \mathcal{O}_{\pi^{-1}(U_i)}$ which glue to a map $\sigma: \pi^*E^{\vee} \ra \mathcal{O}_{V(E)}$ whose zero locus is exactly the zero section $X$ in $V(E)$. From this we get a Koszul resolution of the zero section.
\begin{gather*}
 0 \ra \bigwedge^{\rk E} \pi^*E^{\vee} \ra \cdots \ra \pi^*E^{\vee} \overset{\sigma}{\ra} \mathcal{O}_{V(E)} \ra \mathcal{O}_X \ra 0
 \label{eqn:koszulres}
\end{gather*}
The map $\pi: V(E) \ra X$ is affine so it induces functors $\check{C}(X) \ra \check{C}(V(E))$. Then given a complex of vector bundles $F$ on $X$ the resolution of the complex $i_* F$ is given tensoring $\pi^* F$ with the Koszul resolution and taking the total complex. We shall abuse notation and denote this complex of locally-free sheaves by $i_* F$.

Let $F$ and $D$ be complexes of locally-free sheaves on $X$. Then, there is a
bi-graded dg-algebra associated to $F$ and $D$. Namely,
\begin{equation*}
 V(F,D) = \bigoplus_{l \in \Z} \Hom_X(F,D \otimes \bigwedge^l E)[l]
\end{equation*}
Here the differential respects the grading by exterior powers of $V$. From $V(F,D)$ there is natural map to the endomorphisms of the complex in equation \ref{eqn:koszulres}. We map $\phi \in \Hom_X(F,D \otimes \bigwedge^l E)[l]$ to $\pi^*\phi$ in $\Hom_{V(E)}(\pi^*F,\pi^*D[l])$. Let us name this map $\pi^*$.

\begin{prop}
 $\pi^*$ is a quasi-isomorphism.
 \label{prop:sscompute}
\end{prop}

\proof Note that $\pi^*$ preserves the bi-gradings and thus descends to a morphism
of the spectral sequences associated to the double complex. The $E^2$-pages for each
are given by $E^2_{pq} = H^p(\Hom_X(F,D \otimes \bigwedge^q E)$ and the map induces
an isomorphism. \qed

Thus, we can replace $\Hom_{V(E)}(i_*F,i_*D)$ by $V(F,D)$ and we will. There is one important subtlety - if we consider $V(F,F)$ with the algebra structure induced by $X$, we do not have a morphism of dg-algebras. Thus, we must pullback the algebra structure of $\Hom_{V(E)}(i_*F,i_*F)$ to $V(F,F)$. This is described as follows. Each $\phi \in \Hom_X(F,F \otimes \bigwedge^l E)$ induces $\phi \otimes id \in \Hom_X(F \otimes \bigwedge^k E,F \otimes \bigwedge^{k+l} E)$. This allows us to compose $\psi \circ \phi$ for $\psi \in \Hom_X(F,F \otimes \bigwedge^k E)$. With respect to this algebra structure we have a restatement of the previous proposition.

\begin{cor}
 $\pi^*$ is a quasi-isomorphism of $V(F,F)$ and $\Hom_V(i_*F,i_*F)$ as dg-algebras.
\end{cor}

\begin{rmk}
 This was also observed in the case of line bundles in \cite{Seg07}.
\end{rmk}

We can restate this construction purely in categorical. Let $\mathcal{C}$ be a category with an auto-functor $L: \mathcal{C} \ra \mathcal{C}$. Then, we can construct the trivial extension category of $\mathcal{C}$ by $L$, denoted by $\mathcal{C} \oplus L$. The objects $\mathcal{C} \oplus L$ are the same as the objects of $\mathcal{C}$. The morphism space between objects $A,B$ is defined to be
\begin{gather*}
 \Hom_{\mathcal{C} \oplus L}(A,B) = \Hom_{\mathcal{C}}(A,B) \oplus \Hom_{\mathcal{C}}(A,L(B))
\end{gather*}
With respect to this decomposition, write a morphism $\phi$ as $\phi_1 \oplus \phi_2$. $\phi$ composed with $\psi: B \ra C$ is defined to be $\psi_1 \circ \phi_1 \oplus L(\psi_1) \circ \phi_2 + \psi_2 \circ \phi_1$.

Trivial extensions preserve properties of $\mathcal{C}$ and $L$. By this, we mean that if $\mathcal{C}$ has extra structure and $L$ respects the extra structure, then $\mathcal{C} \oplus L$ also possesses this extra structure. For example, if $\mathcal{C}$ is a dg-category and $L$ is a dg-functor, then $\mathcal{C} \oplus L$ is also a dg-category where the differential is just the direct sum of the differentials on each component of the Hom space.

We have a quasi-equivalence between the image of $i_*$ in $\check{C}(V(E))$ and the trivial extension dg-category $\check{C}(X) \oplus \left(\otimes_{\mathcal{O}_X} \bigwedge^l E [l] \right)$. To ease typographical stress, we shall set $T(E) = \left(\otimes_{\mathcal{O}_X} \bigwedge^l E [l] \right)$. This extends to an equivalence of the appropriate derived categories. 

\begin{lem}
 $i_*$ induces an equivalence between $D(\check{C}(X) \oplus T(E))$ and $D^b_X(\Coh(V(E)))$.
\end{lem}

\proof From the calculation in proposition \ref{prop:sscompute} we see that $i_*$ is full and faithful. Using corollary \ref{cor:generationcor}, we know that the image $i_*$ strongly generates $D_X^b(\Coh(V(E)))$. Since $X$ is smooth, this gives all bounded compelexes of locally-free sheaves. \qed

Now let us simplify matters a little by recalling a result from \cite{BvB02}.

\begin{thm}
 Let $Y$ be a quasi-compact and separated scheme with enough locally-free coherent sheaves. Then, there exists a bounded complex of locally-free coherent sheaves $G$ which generates $D(\Qcoh(Y))$. In particular, if $Y$ is smooth, then $G$ strongly generates $D^b(\Coh(Y))$.
\end{thm}

This means that we only need to pay attention to the dg-algebra $\Hom_{\check{C}(Y)}(G,G)$ (where we have extended the category $\check{C}(Y)$ to include bounded complexes of locally-free coherent sheaves in the obvious way). But, with our assumptions on $X$ (projective over a field $k$) we can assume that $G$ is actually a single locally-free coherent sheaf. The proof in the case does not require the same reduction argument as in \cite{BvB02} and \cite{Nee96}. We give a proof in the appendix.

We can now reduce to a single dg-algebra and modules over it. Let us denote the dg-algebra $V(G,G)$ by $V(X,G)$. From a result of Keller \cite{Kel94}, we have the following.

\begin{prop}
 Since $X$ is smooth, $D^b_X(\Coh(V(E)))$ is triangle equivalent to $D(V(X,G))$.
\label{cor:maincorollary}
\end{prop}

Let us step back for a moment. Take $A$ to be a dg-algebra and $M$ a dg-module over $A$. We can rephrase the construction of the trivial extension categories in more down to earth terms.

\begin{defn}
 Given such an $A$ and $M$ as above we can form the trivial extension dg-algebra of $A$ by $M$ denoted by $A(M)$. It is a dg-algebra over $k$. As a vector space it is just $A \oplus M$. The differential is $d_A \oplus d_M$ and composition is given by $(a,m)\cdot(a'm') = (a \cdot a', a \cdot m' + m \cdot a' + m \cdot m')$.
\end{defn}

\begin{eg}
 A simple of example is given by taking any dg-bi-module $M$ over $A$ and giving it the zero algebra structure, i.e. the composition $M \otimes_A M \ra M$ is just the zero map. Then $A(M)$ is the trivial infinitesimal extension of $A$ by $M$. If $M$ naturally has an algebra structure over $A$, we get a little more structure.
\end{eg}

\begin{rmk}
 It is obvious but worth remarking that we require a bi-module for this construction. Obviously one can promote any left or right dg-module to a bi-module by bestowing it with the trivial action on the right or left, respectively.
\end{rmk}

In general, to determine the structure of the algebra $A(M)$ one needs to identify the bi-module $M$. To explicitly give the algebra, one needs to know how $A$ acts on $M$. There is one special situation where one can determine up to quasi-isomorphism the dg-algebra $A(M)$ solely from the dg-algebra $A$. That is when $M = \Hom_k(A,k) = A^*$. On $A(A^*)$ there is non-degenerate pairing given by $\left< (a_1,b_1^*), (a_2,b_2^*) \right> = \tr(a_1a_2,a_1b_2^*+b_1^*a_2)$ where $\tr: A(A^*) \ra k$ takes $(a,b^*)$ to $(-1)^{|a||b|}b^*(a)$. The pairing satisfies
\begin{equation*}
 \left< (a_1,b_1^*)\cdot(a_2,b^*_2), (a_3,b^*_3) \right> = (-1)^{n_1(n_2+n_3)}\left (a_2,b_2^*)\cdot(a_3,b_3^*), (a_1,b_1^*) \right>
\end{equation*}
where $|a_i|=|b_i|=n_i$ and 
\begin{equation*}
 \left< d(a_1,b^*_1),(a_2,b_2^*) \right> = (-1)^{n_1n_2+1}\left< d(a_2,b_2^*),(a_1,b_1^*) \right>
\end{equation*}

\begin{lem}
 If $A$ is finite-dimensional, then any dg-algebra $B$ with an injective map $A \hookrightarrow B$, a non-degenerate cyclically symmetric pairing as above, and a splitting $B \cong A \oplus I$ with $I^2 = 0$ is isomorphic to $A(A^*)$.
\end{lem}

\proof We have $A^* \cong I$ as chain complexes and we just need to know the values of $a\cdot i$ and $i \cdot a$ for $a \in A$ and $i \in I$. These can found by using the non-degenerate cyclically symmetric pairing. \qed

Unfortunately, we cannot apply this lemma directly to the case of $E = \omega_X$ for two reasons. One, our dg-algebras are not finite dimensional, and two, with our choice of dg-model for $D^b_X(\Coh(V(E)))$ we do not have a non-degenerate pairing, only a pairing which is non-degenerate on cohomology. We can reduce to the cohomology which is finite dimensional over $k$ but we must pass from a dg-algebra to an $\Ainf$-algebra. Through this passage, we can show the following result.

\begin{prop}
 Let $X$ be a smooth projective scheme over a field $k$. Let $G$ be a locally-free coherent sheaf strongly generating $D^b(\Coh(X))$ and $B_X$ the endomorphisms of $G$ in $\check{C}(X)$. Then, $V(X,G)$ is quasi-isomorphic to $B_X(B_X^*[-\dim X - 1])$.  \label{prop:intermediate}
\end{prop}

\begin{cor}
 $D_X^b(\Coh(V(E))$ is triangle equivalent to $D(B_X(B_X^*[-\dim X - 1]))$. \label{cor:mainresult}
\end{cor}

\begin{rmk}
 These results were discussed in the case when the generator chosen for $X$ has no higher cohomology in \cite{Seg07}.
\end{rmk}

\section{A-infinity algebras}

Let $k$ be a field.

\begin{defn}
 Let $V$ be a graded $k$-module. Then bar coalgebra $B(V)$ is the $k$-module $\oplus_{l \geq 0} (sV)^{\otimes l}$ where $(sV)^i = V^{i+1}$. An element $sv_1 \otimes \cdots \otimes sv_k$ will be denoted by $[v_1 | \cdots | v_k]$. $B(V)$ is equipped with the coalgebra structure $\Delta$.
\begin{displaymath}
 \Delta [v_1 | \cdots | v_k] = \sum_{i=0}^k [v_1|\cdots|v_i] \bigotimes [v_{i+1}|\cdots|v_k]
\end{displaymath}
\end{defn}

\begin{defn}
 Given a coalgebra $(C,\Delta)$ a coderivation $d: C \ra C$ is a degree one $k$-linear map for which the following diagram commutes
\begin{center}
\leavevmode
\begin{xy}
 (-15,10)*+{C}="a"; (15,10)*+{C \otimes C}="b"; (-15,-10)*+{C}="c"; (15,-10)*+{C \otimes C}="d"; {\ar@{->}^{\Delta} "a";"c"}; {\ar@{->}^{d} "a";"b"}; {\ar@{->}^{\Delta} "b";"d"}; {\ar@{->}^{d \otimes 1 + 1 \otimes d} "c";"d"}
\end{xy}
\end{center}
$d$ is a codifferential if $d^2=0$.
\end{defn}

The following is well known.

\begin{lem}
 There is bijection between coderivations $B(V)$ and $\bigoplus_l \Hom_k((sV)^{\otimes l},sV)$ given by sending $d$ to $d_l: (sV)^{\otimes l} \ra sV$.
\end{lem}

\begin{defn}
 A curved $\Ainf$-algebra is a $k$-module $A$ with a codifferential $b$ on $B(A)$. If $b_0 = 0$, we shall simply call $A$ an $\Ainf$-algebra.
\end{defn}

\begin{cor}
 The data of an $\Ainf$-algebra on $A$ is equivalent to a collection of degree one maps $b_k: (sA)^{\otimes k} \ra (sA)$ $k > 0$ satisfying
\begin{displaymath}
 \sum_{i+j+k=n,i,k\geq 0,j>0} b_{n-j+1}(\id^{\otimes i} \otimes b_j \otimes \id^{\otimes k})
\end{displaymath}
\end{cor}

\begin{rmk}
 One can also translate this into degree $2-k$ maps $m_k: A^{\otimes k} \ra A$ where $m_k = s^{-1} \circ b_k \circ s^k$. These satisfy
\begin{displaymath}
 \sum_{i+j+k=n,i,k\geq 0,j>0} (-1)^{j+k(n-k-j)}m_{n-j+1}(\id^{\otimes i} \otimes m_j \otimes \id^{\otimes k})
\end{displaymath}
In particular if set $m_k = 0$ for $k>2$. We get 
\begin{displaymath}
 \begin{aligned}
  m_1^2 & = 0 \\
  m_2(m_1,\id) + m_2(\id,m_1) & = 0 \\
  m_2(m_2,\id) - m_2(\id,m_2) & = 0
 \end{aligned}
\end{displaymath}
These are just the requirements for a dg-algebra. In particular, we can convert all the examples from the previous section into this language.
\end{rmk}

\begin{defn}
 A morphism $F: A \ra A'$ of $\Ainf$-algebras is degree zero coalgebra map $F: B(A) \ra B(A')$ that commutes with the codifferentials. $F$ is a called quasi-isomorphism if $F_1$ is a quasi-isomorphism.
\end{defn}

As before, $F$ is determined by the restriction $F_k: (sA)^{\otimes k} \ra sA'$. To recover $F$ we set $F = \sum F_{i_1} \otimes \cdots \otimes F_{i_n}$. For $F$ to commute with differential we need
\begin{displaymath}
 \sum_{i+j+k=n,i,k\geq 0,j>0} F_{n-k+1}(\id_{sA}^{\otimes i} \otimes b^A_j \otimes \id_{sA}^{\otimes k}) = \sum_{i_1+\cdots+i_r=n} b_r(F_{i_1} \otimes \cdots \otimes F_{i_r})
\end{displaymath}

\begin{defn}
 An $\Ainf$-algebra $A$ is called minimal if $b_1=0$.
\end{defn}

$\Ainf$-algebras may seem a bit unwieldy, but there are advantages. One of which is the following result.

\begin{thm}
 Given an $\Ainf$-algebra $A$ and choice of decomposition $A = H(A) \oplus B \oplus D$ where $b_1:D \ra B$ is an isomorphism. There is an $\Ainf$ algebra structure on $H(A)$ and $\Ainf$-quasi-isomorphisms $p: A \ra H(A)$ and $i: A \ra H(A)$ so that the $p_1$ and $i_1$ are the projection and inclusion determined by the decomposition.
\end{thm}

This result is originally due to Kadeishvili \cite{Kad80}. The statement above is taken from \cite{Mar06}. The $\Ainf$-compositions on $H(A)$, the $p_j$, and the $i_j$ are all constructed from iterative applications of $(b_1|_D)^{-1}=h$, the $b_j$ on $A$, $p_1$, and $i_1$.

Let $(C,\Delta)$ be a coalgebra. Then a left comodule $(N,\Delta_N)$ for $C$ is given by $k$-module $N$ and a degree zero $k$-linear map $\Delta_N: N \ra C \otimes N$ so that the following diagram commutes
\begin{center}
\leavevmode
\begin{xy}
 (-15,10)*+{M}="a"; (15,10)*+{C \otimes M}="b"; (-15,-10)*+{C \otimes C}="c"; (15,-10)*+{C \otimes C \otimes M}="d"; {\ar@{->}^{\Delta} "a";"c"}; {\ar@{->}^{\Delta_M} "a";"b"}; {\ar@{->}^{1 \otimes \Delta_M} "b";"d"}; {\ar@{->}^{\Delta \otimes 1} "c";"d"}
\end{xy}
\end{center}
A right comodule is defined similarly. A bicomodule $N$ over $C$ is given by a $k$-module $N$ which is a left comodule under $\Delta_N$ and a right comodule under $\Delta^N$ so that the following diagram commutes
\begin{center}
\leavevmode
\begin{xy}
 (-15,10)*+{N}="a"; (15,10)*+{C \otimes N}="b"; (-15,-10)*+{N \otimes C}="c"; (15,-10)*+{C \otimes N \otimes C}="d"; {\ar@{->}^{\Delta^N} "a";"c"}; {\ar@{->}^{\Delta_N} "a";"b"}; {\ar@{->}^{1 \otimes \Delta^N} "b";"d"}; {\ar@{->}^{\Delta_N \otimes 1} "c";"d"}
\end{xy}
\end{center}

Let $M$ be a $k$-module. We can form a left comodule over $B(A)$, $B(A,M) = B(A) \otimes sM$. The left coaction of $B(A)$ on $B(A,M)$ is given by
\begin{displaymath}
 \Delta_M([a_1|\cdots|a_k|m] = \sum_{i=0}^k [a_1|\cdots|a_i] \bigotimes [a_{i+1}|\cdots|a_k|m]
\end{displaymath}
One can define a right comodule $B(M,A)$ similarly. The coaction here will be denoted by $\Delta^M$.

\begin{defn}
 A left $\Ainf$-module $M$ over an $\Ainf$-algebra $A$ is a $k$-module $M$ equipped with a degree one $k$-linear map $b_M: B(A,M) \ra B(A,M)$ with $b_M^2 = 0$ and the following diagram commuting
\begin{center}
\leavevmode
\begin{xy}
 (-20,10)*+{B(A,M)}="a"; (20,10)*+{B(A,M)}="b"; (-20,-10)*+{B(A) \otimes B(A,M)}="c"; (20,-10)*+{B(A) \otimes B(A,M)}="d"; {\ar@{->}^{\Delta_M} "a";"c"}; {\ar@{->}^{b_M} "a";"b"}; {\ar@{->}^{\Delta_M} "b";"d"}; {\ar@{->}^{b \otimes 1 + 1 \otimes b_M} "c";"d"}
\end{xy}
\end{center}
One can define a right $\Ainf$-module analogously. A morphism $F: M \ra N$ of left $\Ainf$-modules is a $k$-linear degree zero map $F: B(A,M) \ra B(A,N)$ so that $F \circ d_M = d_N \circ F$ and $\Delta^N \circ F = (1 \otimes F) \circ \Delta^M$. A morphism of right $\Ainf$-modules is defined similarly.
\end{defn}

We can also make a $B(A)$-bicomodule from a $k$-module $M$. Let $B(A,M,A) = B(A) \otimes sM \otimes B(A)$. We get a left coaction
\begin{displaymath}
 \Delta_M [a_1|\cdots|a_k|m|a_{k+1}|\cdots|a_{k+l}] = \sum_{i=0}^k [a_1|\cdots|a_i] \bigotimes [a_i|\cdots|a_k|m|a_{k+1}|\cdots|a_{k+1}]
\end{displaymath}
and a right coaction
\begin{displaymath}
 \Delta^M [a_1|\cdots|a_k|m|a_{k+1}|\cdots|a_{k+l}] = \sum_{i=0}^l [a_1|\cdots|a_k|m|a_{k+1}|\cdots|a_{k+i}] \bigotimes [a_{k+i+1}|\cdots|a_{k+l}]
\end{displaymath}

\begin{defn}
 An $\Ainf$-bimodule $M$ over an $\Ainf$-algebra $A$ consists of a $k$-module $M$ with and a degree one $k$-linear $b_M^M: B(A,M,A) \ra B(A,M,A)$ so that $(b_M^M)^2 = 0$ and the following diagrams commute
\begin{center}
\leavevmode
\begin{xy}
 (-25,10)*+{B(A,M,A)}="a"; (-25,-10)*+{B(A,M,A)}="c"; (25,10)*+{B(A)\otimes B(A,M,A)}="b"; (25,-10)*+{B(A)\otimes B(A,M,A)}="d"; {\ar@{->}^{b_M^M} "a";"c"}; {\ar@{->}^{\Delta_M} "a";"b"}; {\ar@{->}^{(b\otimes 1 + 1 \otimes b_M^M)} "b";"d"}; {\ar@{->}^{\Delta_M} "c";"d"}
\end{xy}
\end{center}
\begin{center}
\leavevmode
\begin{xy}
 (-25,10)*+{B(A,M,A)}="a"; (-25,-10)*+{B(A,M,A)}="c"; (25,10)*+{B(A,M,A)\otimes B(A)}="b"; (25,-10)*+{B(A,M,A)\otimes B(A)}="d"; {\ar@{->}^{b_M^M} "a";"c"}; {\ar@{->}^{\Delta^M} "a";"b"}; {\ar@{->}^{(b_M^M\otimes 1 + 1 \otimes b)} "b";"d"}; {\ar@{->}^{\Delta^M} "c";"d"}
\end{xy}
\end{center}
\end{defn}

\begin{rmk}
 Since $\Ainf$-bi-modules will be of interest, we shall labor this point for a moment longer. $b^M_M: B(A,M,A) \ra B(A,M,A)$ is uniquely determined by the maps $b_{k,l}: (sA)^{\otimes k} \otimes (sM) \otimes (sA)^{\otimes l} \ra sM$. These maps must satisfy
\begin{displaymath}
 \sum_{i + j + k =n,i,j \geq 0, k>1} b_{n-k+1,m}(\id_{sA}^{\otimes i} \otimes b_k \otimes \id_{sA}^{\otimes j} \otimes \id_{sM} \otimes \id_{sA}^{\otimes n}) +
\end{displaymath}
\begin{displaymath}
 \sum_{s+i=n,h+j=m,\ i,j,h,s\geq 0} b_{n-i,m-j}(\id_{sA}^{\otimes l} \otimes b_{i,j} \otimes \id_{sA}^{\otimes h}) +
\end{displaymath}
\begin{displaymath}
 \sum_{s+r+t=m,s,t\geq 0,r>0} b_{n,m-t+1}(\id_{sA}^{\otimes n} \otimes \id_{sM} \otimes \id_{sA}^{\otimes r} \otimes b_t \otimes \id_{sA}^{\otimes s}) = 0
 \label{rmk:relations}
\end{displaymath}
If we set $b_{i,j}=0$ for $i+j>1$ and $b_k=0$ for $k>1$, we get the following relations.
\begin{displaymath}
 \begin{aligned}
  b_{0,0}^2 & = 0 \\
  b_{1,0}(b_1 \otimes \id) + b_1(b_{1,0}) & = 0 \\
  b_{0,1}(\id \otimes b_1) + b_1(b_{0,1}) & = 0 \\
  b_{0,1}(b_{1,0} \otimes \id) + b_{1,0}(\id \otimes b_{0,1}) & = 0 \\
  b_{1,0}(\id \otimes b_{1,0}) + b_{1,0}(b_2 \otimes \id) & = 0 \\
  b_{0,1}(b_{0,1} \otimes \id) + b_{0,1}(\id \otimes b_2) & = 0
 \end{aligned}
\end{displaymath}
After shifting the gradings back, we get the relations satisfied by a dg-bi-module over a dg-algebra.
\end{rmk}

\begin{eg}
\begin{enumerate}
 \item Given an $\Ainf$-algebra $A$, we can define an $\Ainf$-bimodule structure on $A$ by setting $b_{k,l} = b_{k+l+1}$.
 \item From an $\Ainf$-bi-module $M$, we can construct an $\Ainf$-bimodule structure on $M^* = \Hom_k(M,k)$ as follows. Consider $B(A,M^*,A) \otimes sM$ and let $\sigma$ denote the $k$-module isomorphism $sM^* \otimes B(A,M,A)$ given by shifting the first copy of $B(A)$ to the end. Let $\tr: M^* \otimes M \ra k$ denote the natural pairing. We can define $b^{M^*}_{M^*}$ implicitly by setting $\tr(b^*_{k,l} \otimes \id) = \tr(\id \otimes b_{k,l} \circ \sigma)$. One can quickly check that $b^*_{k,l}$ satisfy the relations \ref{rmk:relations}. Explicitly, we have
\begin{displaymath}
 b^*_{k,l}([a_1|\cdots|a_k|m^*|a_{k+1}|\cdots|a_{k+l}])(m') = (-1)^{\bowtie}m^*(b_{l,k}([a_{k+1}|\cdots|a_{k+l}|m'|a_1|\cdots|a_k])
\end{displaymath}
where $\displaystyle{\bowtie = (\sum_{r=1}^k|a_r|)(|m^*|+\sum_{t=1}^l|a_{k+t}|)+|m^*|+\sum_{j=1}^k|a_j|(\sum_{i=j+1}^k|a_i|)}$
\end{enumerate}
\end{eg}

Let $M$ be an $k$-module. Let us set 
\begin{displaymath}
\underbrace{B(A) \otimes M \otimes B(A) \otimes M \otimes \cdots \otimes M \otimes B(A)}_{n \text{ copies of } M} =  B_n(A,M,A) 
\end{displaymath}
Then,
\begin{displaymath}
 B(A\oplus M) = \bigoplus_{n \geq 0} B_n(A,M,A)
\end{displaymath}
Note that $B_1(A,M,A) = B(A,M,A)$. 

\begin{defn}
 An $\Ainf$-algebra $A'$ over an $\Ainf$-algebra $A$ is codifferential $b'$ on $B(A \oplus A')$ whose restriction to $B(A)$ is $b$.
\end{defn}

The following is clear.

\begin{lem}
 Given an $\Ainf$-bimodule $M$ over an $\Ainf$-algebra, we can form the trivial extension $\Ainf$-algebra $A(M)$ which as $k$-module is $A \oplus M$ and which has as codifferential $b_{A(M)} = b + b^M_M$ on $B(A \oplus M)$. This is an $\Ainf$-algebra over $A$.
\end{lem}

Let $\Mod-A$ denote the dg-category of right $\Ainf$-modules over an $\Ainf$-algebra $A$. Morphisms from $M$ to $N$ are all $k$-linear maps from $B(A,M)$ to $B(A,N)$ commuting with comodule structure. The differential is given by commuting a given morphism with the differentials. This dg-category is pre-triangulated with triangles coming from cones in the dg-category. Let $D(A)$ denote the smallest triangulated subcategory of $H^0(\Mod-A)$ containing $A$ and closed under idempotent splittings. The following result is useful.

\begin{thm}
 If $A$ and $B$ are quasi-isomorphic $\Ainf$-algebras, then $D(A) \cong D(B)$.
\end{thm}

For a proof, see \cite{SeiDR}.

Previously, we had defined the bounded derived category of special dg-categories, ones where all quasi-isomorphisms factor through homotopy equivalences. The category of right dg-modules sits inside $\Mod-A$.

\begin{prop}
 Given a dg-algebra $A$ as above, then the two definitions of the derived category of $A$ coincide.
\end{prop}

\proof Both these categories are equivalent to $D(\dgMod-\Omega B A)$ where $\Omega B A$ is the bar-cobar of $A$, see \cite{LH03}. \qed

\begin{defn}
 A cyclic $\Ainf$-algebra of dimension $d$ is an $\Ainf$-algebra
$A$ with a symmetric, cohomologically non-degenerate pairing $(\cdot,\cdot): A \otimes A
\ra k$ of degree $-d$ such that 
\begin{equation*}
 \left(m_n(a_0,\cdots,a_{n-1}),a_n \right) = (-1)^{n+|a_0|(|a_1|+\cdots+|a_n|)}\left( m_n(a_1,\cdots,a_n),a_0 \right)
\end{equation*}
for all $n$.
\end{defn}

\begin{rmk}
 \begin{enumerate}
 \item If we using the gradings on $sA$ instead, the pairing becomes honestly cyclically symmetric.
 \item We can turn the pairing into a A-bimodule map $f: A \ra A^*[-d]$ by setting $f_1(a) = (a,-)$ and $f_n=0$. Nondegeneracy says this is a quasi-isomorphism.
 \end{enumerate}
\end{rmk}

Now we prove our main technical lemma.

\begin{lem}
 Let $A'$ be an $\Ainf$-algebra over $A$ both having finite dimensional cohomology. Assume that $A'$ is cyclic, that $b_{A'}$ respects the auxiliary grading $B(A\oplus A') = \oplus B_n(A,A',A)$, and the pairing $(\cdot,\cdot)$ has auxiliary degree zero with the auxiliary degree of the target field being one. Then, $A'$ is quasi-isomorphic to $A(A^*)$.
\label{lem:technicalresult}
\end{lem}

\proof Appealing to proposition \ref{prop:niceminimalmodel}, we see that it is enough to consider the situation where $A'$ and $A$ are minimal and finite dimensional. The pairing gives an isomorphism $A^* \cong A'$. Cyclicity of $(b_k,\id)$ and isotropy of $(\cdot,\cdot)$ force $b_k$ to vanish on $B_n(A,A',A)$ for $n > 1$ to map $B_i(A,A',A)$ into $B_i(A,A',A)$ for $i=0,1$. They also force $b_{k,l}: A^{\otimes k} \otimes A' \otimes A^{\otimes l} \ra A'$ to coincide under the isomorphism $A^* \cong A'$ with the $\Ainf$-bimodule structure on $A^*$. Thus, we see that $A'$ is isomorphic to $A(A^*)$ in this case and, hence, the result holds in general. \qed

\begin{rmk}
 One can view proposition \ref{lem:technicalresult} as a slightly weaker characterization of the ``dualizing bimodule'' for an $\Ainf$-algebra.
\end{rmk}

\begin{rmk}
 Preservation by $b$ of the decomposition $B(A\oplus A') = \bigoplus B_n(A,A',A)$ is equivalent to the existence of an action by $k^{\times}$ on $B(A \oplus A')$ which sends $b = \oplus b|_{B_n(A,A',A)}$ to $b' = \oplus \lambda^n b|_{B_n(A,A',A)}$. This action can be realised geometrically in the case of $V \ra X$ as the standard action of $k^{\times}$ on a vector bundle. This was pointed out to me by Paul Seidel.
\end{rmk}

\section{Formal non-commutative symplectic geometry}

In this section, we prove Proposition \ref{prop:niceminimalmodel}. Before we dive in, we review the notions of formal non-commutative symplectic geometry and its relation with $\Ainf$-algebras. Better references for this are the original source \cite{Kon93} and \cite{Gin01,HL04}. I refer the reader to one of these for any unproven assertions.

Let $k$ be field for simplicity.

\begin{defn}
 A formal $k$-module $V$ is a $\Z$ (or $\Z/2\Z$)-graded $k$-module arising as the inverse limit of finite-dimensional $k$-modules $V_i$ and equipped with the inverse limit topology.
\end{defn}

\begin{defn}
 A topological basis for $V$ is a choice of basis of $V$ for which any element can be written uniquely as a convergent sum.
\end{defn}

The category of formal $k$-modules will be denoted $f-k-\Mod$. Morphisms in this category are continuous $k$-module morphisms. Given two formal $k$-module $V$ and $W$, we form the completed tensor product $V \hat{\otimes} W := \underset{\leftarrow}{\lim} \ V_i \otimes_k W_j$. 

\begin{prop}
 The functors $k-\Mod \ra f-k-\Mod$, $V \mapsto V^*$ and $f-k-\Mod \ra k-\Mod$, $W \mapsto
W^{\star}$ where $W^{\star}$ is the topological dual, are contravariant equivalences of symmetric monoidal categories.
\end{prop}

A formal $k$-algebra $A$ will be an associative, unital algebra object in this category. We wish to treat $A$ as the ring of functions on some un-named type of manifold. Any reader with an acquaintance with differential geometry should recognize the similarity. To wit, we employ the following notation.

\begin{defn}
 A continuous derivation $X: A \ra A$ is called a vector field. A continuous algebra homomorphism $\phi: A \ra B$ is called a smooth map. If it is invertible, we call it a diffeomorphism.
\end{defn}

\begin{eg}
 Let $V$ be a graded $k$-module and consider the tensor algebra on $V$, $\bigoplus_{i \geq 0} V^{\otimes k}$ Then, the dual considered a formal $k$-module can be written as $\hat{T}V^{\vee} = \prod_{i \leq 0} V^{\hat{\otimes} i}$. Dualising the co-algebra structure gives the standard algebra structure on a tensor algebra. This will be our standard example of a formal $k$-algebra. Note that we can view it as formal non-commutative power series the variables $x^s$ for $s \in S$ a basis for $V$.
\end{eg}

Assume our formal algebra $A$ is isomorphic to $\hat{T}V$ for some $V$. Then there are two gradings on $A$. The grading induced by $V$ where $x^1 \otimes \cdots \otimes x^n$ has degree $|x^1|+\cdots+|x^n|$ and the grading induced by the tensor algebra where $x^1 \otimes \cdots \otimes x^n$ has degree $n$. The second grading we shall call the order.

Any vector field $X$ on $\hat{T}V$ can be expanded as $X = \sum_{i \geq 0} X_i$ where $X_i: V \ra V^{\hat{\otimes} i-1}$. We say a vector field vanishes at zero if $X_0 = 0$. A smooth map $\phi: \hat{T}V \ra \hat{T}W$ can be also be expanded $\sum_{i \geq 1} \phi_i$ with $\phi_i: V \ra W^{\hat{\otimes} i}$.

\begin{lem}
 A smooth map $\phi: \hat{T}V \ra \hat{T}W$ is a diffeomorphism if and only if $\phi_1$ is invertible.
\end{lem}

\begin{defn}
 The module of one-forms $\Omega^1(A)$ is defined to be $A \hat{\otimes} A/k$. We write $x \hat{\otimes} y$ as $xdy$. Correspondingly, there is a natural map $d: A \ra \Omega^1(A)$. $\Omega^1(A)$ possesses the structure of an $A$-bimodule with $a \cdot xdy = axdy$ and $xdy \cdot a = xd(ya)-xyda$.
\end{defn}

\begin{prop}
 There is a bijection between derivations of $A$ with values in $M$ and continuous maps $\Omega^1(A)$ to $M$.
\end{prop}

\begin{defn}
 The dg-algebra of tensor forms is $\Omega^*(A)$ the completed tensor algebra over $A$ on $s^{-1}\Omega^1(A)$ where $s$ is the parity shift, $(sW)^i = W^{i+1}$. It possesses the standard algebra structure. $d$ extends to a differential on $\Omega^*(A)$.
\end{defn}

Given a vector field $X$ on $A$, we can define some operations on $\Omega^*(A)$. Let $i_X$ be the derivation of degree $X$ of $\Omega^*(A)$ for $i_X (a) = 0$ and $i_X (da) = X(a)$ and let $L_X$ be the derivation of degree $X$ of $\Omega^*(A)$ for which $L_X(a) = X(a)$ and $L_X(da) = d(L_X(a))$. These are the contraction and Lie derivative.

\begin{defn}
 The de Rham complex of $A$ is $DR^*(A) = \Omega^*(A)/[\Omega^*(A),\Omega^*(A)]$. The differential $d$ and derivations $i_X$ and $L_X$ descend to $DR^*(A)$. 
\end{defn}

The standard relations hold. Here are two.

\begin{lem}
 $L_X = [i_X,d], [L_X,L_Y] = L_{[X,Y]}$
\end{lem}

Choose a topological basis $x^i$ of $A$. Then products of $x^i,dx^i$ form a topological basis for $DR^*(A)$. Let $(DR^*(A))^{[p]}$ denote the subcomplex of $DR^*(A)$ formed from forms with tensor order divisible by $p$.

\begin{lem}
 (Poincar\'{e}) $H^*(DR^*(A),d) = H^*((DR^*(A))^{[p]},d)$ if the characteristic of $k$ is $p$. In particular, if $\charac k = 0$, then $DR^*(A)$ is a resolution of $k$.
\end{lem}

\proof Consider the Euler vector field $E = \sum x^i \del_{x^i}$. The action of $L_E$ on $DR^*(A)$ is diagonalizable with non-negative integral eigenvalues. If the eigenvalue associated to $\nu$ is non-zero and $\nu$ is closed, the Cartan formula for $L_E$ shows that $\nu$ is exact. The zero eigenvalues occur only on $DR^*(A)^{[p]}$. \qed

\begin{defn}
 An $\omega \in DR^2(A)$ is called symplectic if the map $X \ra i_X\omega$ is a bijection and $d\omega = 0$.
\end{defn}

\begin{lem}
 Assume that $A \cong \hat{T}V$ with $V$ finite-dimensional. A two-form $\omega = \sum a_{ij} dx^idx^j$ on $A$ is symplectic if and only if the pairing $V^{\vee} \otimes V^{\vee} \ra k$ given by $(a,b) \mapsto i_{\del_a}i_{\del_b}\omega$ is non-degenerate.
\end{lem}

\proof This is clear. \qed

\begin{cor}
 With the hypotheses as in the previous lemma, we see that $\omega$ is symplectic if and only if $\omega_0$ is symplectic.
\end{cor}

\begin{lem}
 Assume that $A \cong \hat{T}V$. Given a topological basis $x^a$ for $V$, then $x^a$ and $dx^a$ generate $DR^*(A)$.
\end{lem}

\proof Using the commutation relations $d(ab) = dab-(-1)^{|a|}db$ we see that we can reduce the order of the one-form. \qed

If $A \cong \hat{T}V$ we can write any two-form $\omega = \sum_{i \geq 0} \omega_i$ where $\omega_i$ has order $i+2$.

\begin{prop}
 (Darboux) Let $A \cong \hat{T}V$ be a formal $k$-algebra with a symplectic form $\omega$, which is exact. Assume that $V$ is finite dimensional. Then, there exists a diffeomorphism $\phi: A \ra A$ so that $\phi^*\omega = \omega_0$.
\end{prop}

\proof Let us first attempt to solve $L_X \omega_0 = \omega_i$. Note that since $V$ is finite dimensional $\omega_0$ is symplectic if and only if $\omega$ is. Now since $\omega$ is exact so is $\omega_i$ for each $i$. We just need to solve $i_X \omega_0 = \alpha_i$ for some $X$ which we can since $\omega_0$ is symplectic. Let $X^1$ denote the solution to $L_{X^1}\omega_0 = \omega_1$ and consider the diffeomorphism determined by sending $x^i \in V$ to $x^i - V(x^i)$. The order zero term in $\phi^*\omega$ is $\omega_0$ while the order one term is $\omega_1 - L_V\omega_0$ which is zero by construction. Iterating this procedure and noting that the composition of the resulting maps is convergent, we get the desired diffeomorphism $\phi: A \ra A$. \qed

\begin{defn}
 An $\Ainf$-algebra is a formal $k$-algebra $A \cong \hat{T}V$ for some $V$ with a degree one vector field $m: A \ra A$ vanishing at zero and satisfying $[m,m] = 0$.
\end{defn}

By taking the dual of the description of an $\Ainf$-algebra involving the bar complex, we see that we get exactly the above definition.

\begin{defn}
 A symplectic $\Ainf$-algebra is an $\Ainf$-algebra $A$ with a symplectic form $\omega$ so that $L_m \omega = 0$. A homologically-symplectic $\Ainf$-algebra is an $\Ainf$-algebra with closed $m$-constant two-form $\omega$ so that the form induced by $\omega_0$ on the cohomology of $A$ is symplectic.
\end{defn}

\begin{eg}
 Let $B$ be a dg-algebra over a field with map $\tr: B \ra k$ so that the pairing $\left< a,b \right> = \tr(ab)$ is cyclic and non-degenerate on cohomology. Then, the corresponding formal $k$-algebra is a homologically symplectic $\Ainf$-algebra with $\omega$ determined as follows. Express $\tr: B \ra k$ in terms of a basis for $B$ as $\tr(x_i) = t_i$ and express composition as $x_ix_j = \sum_k m^k_{ij} x_k$. Then, $\omega = t_km^k_{ij}dx^idx^j$. That it is $L_m$-constant follows from the proceeding calculation. From the Poincar\'{e} lemma, we see that assumption on the characteristic of $k$ forces $\omega$ to be exact.
\end{eg}

\begin{lem}
 Given a two form $\omega = \sum a_{ij} dx^i dx^j$ of order two, $L_m \omega = 0$ if and only if the corresponding pairing is cyclic.
\end{lem}

\proof  Using the commutation relation, we can write $\omega=\sum_{i \leq j} a_{ij}dx^idx^j$ for some ordering of a topological basis of $V$ corresponding to a basis of the topological dual to $V$. Then the pairing
\begin{equation*}
 \left< x_i,x_j \right> = \begin{cases} a_{ij} & \text{ if $i < j$} \\
                           2 a_{ii} & \text{ if $i = j$} \\
			   (-1)^{|x_i||x_j|}a_{ji} & \text{ if $i > j$}
                          \end{cases}
\end{equation*}
We can write $m x^k = \sum m^k_{i_1\cdots i_l} x^{i_1}\cdots x^{i_l}$. Then, we dualise and use the suspended degree
\begin{equation*}
 \left< m_n(x_{i_0}\cdots x_{i_{n-1}}),x_{i_n} \right> = \left< m^k_{i_0\cdots i_{n-1}}x_k,x_{i_n} \right> = 
\end{equation*}
\begin{equation*}
\sum_{k \leq i_n} m^k_{i_0\cdots i_{n-1}} a_{ki_n} + \sum_{k \geq i_n} (-1)^{(k+1)(i_n+1)}m^k_{i_0 \cdots i_{n-1}} a_{i_n k}
\end{equation*}
The order $n$ term in $L_m \omega$ comes from $L_{m_i}\omega = \sum_{i 
\leq j}a_{ij}m^i_{l_1\cdots l_n}d(x^{l_1}\cdots x^{l_n})dx^j + (-1)^i a_{ij}m^j_{k_1\cdot k_i} dx^id(x^{k_1} \cdots x^{k_i})$.
The term corresponding to $x^{t_0}\cdots dx^{t_{a}} \cdots dx^{t_n}$ carries coefficients
\begin{equation*}
 \sum_{i \leq t_n} (-1)^{|t_1|+\cdots+|t_{a-1}|} a_{it_n}m^i_{t_0\cdots t_{n-1}} +
\end{equation*}
\begin{equation*}
\sum_{i \leq t_a} (-1)^{|t_{a+1}|+\cdots+|t_n|+\sum_{s=a+1}^{n-1}|t_s|(|t_1|+\cdots+\widehat{|t_s|}+\cdots+|t_n|)} a_{it_a}m^i_{t_{a+1}\cdots t_{a-1}} +
\end{equation*}
\begin{equation*}
\sum_{i \geq t_n} (-1)^{i+|t_1|+\cdots+|t_{a-1}|+(|t_n|+1)(|t_0|+\cdots+|t_{n-1}|+1)} a_{t_n i} m^i_{t_0\cdots t_{n-1}} +
\end{equation*}
\begin{equation*}
\sum_{i \geq t_a} (-1)^{\circledS(i)}a_{t_a i} m^i_{t_{a+1}\cdots t_{a-1}}
\end{equation*}
where $\circledS(i) = i+|t_{a+1}|+\cdots+|t_{n-1}|+(|t_a|+1)(|t_1|+\cdots+\widehat{|t_a|}+\cdots+|t_n|+1)+(|t_n|+1)(|t_0|+\cdots+|t_{n-1}|+1)+\sum_{s=a+1}^{n-1}|t_s|(|t_1|+\cdots+\widehat{|t_s|}+\cdots+|t_n|)$. We see this being zero for any $a$ is equivalent to
\begin{equation*}
 \left< m_n(x_{t_0},\cdots,x_{t_{n-1}}),x_{t_n} \right > = (-1)^{|t_0|(|t_1|+\cdots+|t_n|)}\left < m_n(t_1,\cdots,t_n), t_0 \right>
\end{equation*}
using the suspended gradings. If we translate this into the grading on $V$, we get that
\begin{equation*}
 \left< m_n(x_{t_0},\cdots,x_{t_{n-1}}),x_{t_n} \right > = (-1)^{n+|t_0|(|t_1|+\cdots+|t_n|)} \left < m_n(x_{t_1},\cdots,x_{t_n}), x_{t_0} \right>
\end{equation*} \qed

\begin{lem}
 Let $A$ be an $\Ainf$-algebra with a homologically-nondegenerate, $L_m$-constant, exact two-form $\omega$. Write $\omega = \sum_{i \geq 0} \omega_i$. Then, there is a minimal symplectic $\Ainf$-algebra $A'$ quasi-isomorphic to $A$ with $\omega_0|_{H^*}$ as the symplectic form.
\end{lem}

\proof Let $A = TV^*$ and split $V = H \oplus B \oplus D$ and take the morphism coming from the minimal model algorithm, $p:A \ra \hat{T}H^*$, then $p^*\omega$ is non-degenerate on $H^*$ since the constant piece is simply the restriction of $\omega$ to $\hat{T}H^*$ and $p^*\omega$ is $L_{m'}$-constant for the induced $\Ainf$-structure $m'$ on $H^*$ since $p^*L_m = L_{m'}p^*$. Now, apply the Darboux lemma to $A'$ with $p^*\omega$ to get a diffeomorphism $\phi: A' \ra A'$ with $\phi^*p^*\omega = \omega_0|_{H^*}$ and take the vector field given by $\phi \circ m' \circ \phi^{-1}$ as the new $\Ainf$-structure which gives the result. \qed

\begin{rmk}
 This result appears in \cite{Laz06}. Older results deal only with symplectic structures (or at least non-degenerate pairings at the chain level) which one rarely sees. One can turn this into an algorithm similar to \cite{Mar06} modulo the choices of $\alpha_i$ in the Darboux lemma.
\end{rmk}

\begin{cor}
 Let $X$ be a smooth projective Calabi-Yau variety over a field $k$. Then, $D^b(\Coh(X))$ is equivalent to $D(A)$ for a cyclic $\Ainf$-algebra $A$.
\end{cor}

Let us now revisit the previous situation with a new assumption. Assume that we have an auxiliary $\Z$-grading on our $\Ainf$-algebra with $m$ degree zero for this new grading, and $\omega$ degree $l$ with respect to the new induced grading on $DR(A)$. The grading then descends to the cohomology $H(A)$.

\begin{lem}
 There is a minimal $\Ainf$-structure on $H(A)$ which respects the auxiliary grading and for which the $A$ is quasi-isomorphic to $H(A)$ via $\Ainf$-homomorphisms of degree zero.
\end{lem}

\proof Decompose $A$ with respect to the auxiliary grading $A = \bigoplus_{i \in \Z} A_i$ and choose a splitting of $A = H \oplus B \oplus D$ which respects this splitting, i.e. $A_i = H_i \oplus B_i \oplus D_i$ and $d: D_i \ra B_i$ is an isomorphism. Then, the inverse to $d$ also preserves the auxiliary grading. Inspecting the algorithm for constructing the minimal $\Ainf$-structure $H(A)$ and the quasi-isomorphisms, we see that they are formed from iterated compositions of $h, m_i, \pi, i$ where $\pi: A \ra H$ and $i: H \ra A$ hence they are of degree zero. \qed

As a corollary of this, we see that the induced symplectic form, which shall also be denoted by $\omega$, is still of degree $l$. In particular, each term in the order expansion $\omega = \sum \omega_i$ is of degree $l$. Now, our minimal model is finite dimensional over $k$.

\begin{lem}
 The vector fields constructed in the Darboux lemma can be chosen to have degree zero.
\end{lem}

\proof Assume we are in the situation of induction step in the proof of Darboux lemma. Namely, we have a degree $l$ symplectic form $\omega$ on $A$ with an expansion $\omega = \omega_0 + \omega_i + \omega_{i+1} + \cdots$ by order. We seek a vector field $X_i$ so that $L_{X_i} \omega_0 = \omega_i$. We reduced to solving $i_{X_i} \omega_0 = \alpha_i$ where $d\alpha_i = \omega_i$. Decomposing $\alpha_i$ with respect to the auxiliary grading and applying $d$ we see that components with degree not equal to $l$ do not contribute. Hence, we can assume that $\alpha_i$ has degree $l$. One can easily check, that given a vector field $Y$ of degree $t$ and order $i-1$, contracting with $\omega_0$ produces a one-form of degree $t$ and order $i$. Since contraction with $\omega_0$ furnishes an isomorphism of the space of vector fields and one-forms, our $X_i$ can be taken to be of degree zero. \qed

\begin{cor}
 The diffeomorphism constructed in the Darboux lemma can be taken to have degree zero with respect to the auxiliary grading.
\end{cor}

Let us collect these results into a useful conclusion.

\begin{prop}
 Let $A$ be an homologically-symplectic $\Ainf$-algebra with an auxiliary $\Z$-grading so that $m$ is of auxiliary degree zero and the homologically-symplectic form $\omega$ is degree $l$. Then, we can find a symplectic minimal model of $A$ with an auxiliary $\Z$-grading for which the homological vector field is of auxiliary degree zero and the symplectic form is constant and of degree $l$. Moreover, the quasi-isomorphisms are of auxiliary degree zero.
\label{prop:niceminimalmodel}
\end{prop}

Combining this with result with the our lemma \ref{lem:technicalresult} proves proposition \ref{prop:intermediate} and hence the main result \ref{cor:mainresult}.

\section{Closing Remarks}

A useful perspective to take on the results of this paper is furnished by the idea of homological mirror symmetry, \cite{Kon95}. Here we shall frame it rather loosely and optimistically. Given a category arising from algebraic geometry, there is another construction of this category in the realm of symplectic geometry. For example, take a smooth projective Calabi-Yau $d$-fold over $\C$. From this we can construct $D^b(\Coh(X))$. Then, for general enough $X$, homological mirror symmetry predicts that there is another smooth projective Calabi-Yau $d$-fold $Y$ and a category made from it using symplectic geometry, the Fukaya category $\Fuk(Y)$. It is in fact an $\Ainf$-category whose objects are well-mannered and properly groomed Lagrangian submanifolds and whose morphisms between such $L$ and $L'$, when these Lagrangians intersect transversely, are the intersection points. The $\Ainf$-compositions $m_n: \Hom(L_1,L_2) \otimes \cdots \otimes \Hom(L_n,L_{n+1}) \ra \Hom(L_1,L_{n+1})$ come from counting pseudo-holomorphic polygons with appropriate boundary conditions. After we derive $\Fuk(Y)$, we should obtain a category triangle equivalent to $D^b(\Coh(X))$. And $D^b(\Coh(Y))$ should also be equivalent to $D(\Fuk(X))$.

The conjecture extends beyond Calabi-Yau varieties however we lose a bit of the symmetry. A case with non-trivial canonical bundle is that of projective space itself $\P^n_{\C}$. We have two categories of interest $D^b(\Coh(\P^n_{\C}))$ and $D(\Fuk(\P^n_{\C})$. The mirror of this space is not another variety but a variety with a function on it, also called a Landau-Ginzburg model. For $\P^n_{\C}$, the mirror $W: (\C^{\times})^n \ra \C$ is given by
\begin{displaymath}
 W(z_1,\cdots,z_n) = z_1 + \cdots + z_n + \frac{1}{z_1\cdots z_n}
\end{displaymath}
Associated to a Landau-Ginzburg model $W: Y \ra \C$, there are also two categories. The first is analogous to the Fukaya category. Let us assume that $0$ is regular value, then we take as objects Lagrangians with boundary along $W^{-1}(0)$ whose image under $W$ near $0$ is a curve that ends at $0$. The morphisms are again intersection points, subject to some ordering prescriptions near $0$, and compositions arise in the same manner as before. This category symplectically measures the singularities of $W$. In the case of the mirror to $\P^n_{\C}$, we have the Lagrangian thimbles associated to the non-degenerate critical points of $W$ as the representative object. For the second category, we measure the singularities of $W$ algebraically. We consider the fibers $W_{\lambda}$ for $\lambda \in \C$ and take the Verdier quotient of $D^b(\Coh(W_{\lambda})$ by the smallest triangulated category containing all locally-free sheaves and denote it by $D_{Sing}(W_{\lambda})$. Then, the category of singularities for $W$, $D_{Sing}(W)$ is then
\begin{displaymath}
 D_{Sing}(W) = \prod_{\lambda \in \C} D_{Sing}(W_{\lambda})
\end{displaymath}
Homological mirror symmetry predicts that the we again exchange the algebro-geometric categories and symplecto-geometric categories under the duality. In each of these case, one hopes that there are a few central geometrically-interesting objects on each side whose endomorphism $\Ainf$-algebras are quasi-isomorphic.

Now, the evidence for homological mirror symmetry is not overwhelming but it is quite convincing. In the case of Calabi-Yau varieties, little is proven, outside the case of abelian varieties and K3 surfaces. For the Fano or near Fano case, the correspondence between Lagrangians on the variety and the category of singularities of its Landau-Ginzburg mirror is still wide open, but a lot is known about the correspondence between coherent sheaves on the varieties and Lagrangians in Landau-Ginzburg models. Surfaces and toric varieties in general have been thoroughly studied. In this case, there are a few simplifying properties of the geometry of a Landau-Ginzburg model that allow one to bypass a lot of the technical machinery that is needed in the general case.

How does this paper fit in? We seek to build upon on the solid ground of the known cases for surfaces and toric varieties. Given a variety $X$, we extract a Calabi-Yau in a few ways. The most common is to take a divisor in the anti-canonical class if $-\omega_X$ is effective. One could also take affine open set $X$. The third construction is the one appearing in this paper - taking the total space of the canonical bundle. The results here characterize the passage on the level of $\Ainf$-algebra/triangulated categories from $X$ to $\omega_X$ purely in terms of ($\Ainf$-)algebra. One simply takes the trivial nilpotent extension by the dual bi-module. If homological mirror symmetry is true, we should hope that there is a symplecto-geometric method for achieving this algebraic transition. Moreover, believing homological mirror symmetry, we already know the final answer.

\appendix

\section{Generators for projective schemes over a field}

The results of this section are due to Kontsevich although the author learned
about them (as he learned most of his maths, adviser aside) from a paper of
Seidel \cite{Sei03}. For notions of generation see \cite{BvB02} or \cite{Rou03}.

Let $X$ be a projective scheme over a field $k$. Choose some embedding $i:X \ra \P^N_k$. We consider the functor $Li^*: D(\Qcoh(\P^N)) \ra D(\Qcoh(X))$. It commutes with coproducts because it has right adjoint.

\begin{prop}
 The collection $\{\mathcal{O}_X,\ldots,\mathcal{O}_X(N)\}$ generates $D(\Qcoh(X))$.
\end{prop}

\proof Using the Koszul resolution on $\P^N_k$, we can get any tensor power of the twisting sheaf from iterated cones over objects in the collection $\lbrace \mathcal{O}_{\P^N_k}, \ldots, \mathcal{O}_{\P^N_k}(N) \rbrace$.  Pulling back via $Li^*$ shows that we can get any $\mathcal{O}_X(j)$. Using Serre's theorem, we can construct a locally-free quasi-coherent resolution of any quasi-coherent sheaf on $X$ from the set $\lbrace \mathcal{O}_X(j) \rbrace_{j \in \Z}$. Since $Li^*$ commutes with coproducts, the image of $Li^*$ is closed under direct sums. Using homotopy limits we can get any bounded above complex from its truncations. Therefore, we can get all quasi-coherent sheaves starting from $\{\mathcal{O}_X,\ldots,\mathcal{O}_X(N)\}$. The smallest triangulated subcategory containing all quasi-coherent sheaves and closed under direct sums and shifts is all of $D(\Qcoh(X))$. \qed

Let $G = \bigoplus_{j=0}^N \mathcal{O}_X(j)$. Then, $G$ generates $D(\Qcoh(X))$ since we can split idempotents, see \cite{BN93}.

\begin{cor}
 If $X$ is a projective scheme over $k$, then there exists a vector bundle $G$ that generates $D(\Qcoh(X))$ up to idempotent splittings.
\end{cor}

\begin{cor}
 If $G'$ is any generator (possibly up to idempotent splittings) for $D(\Qcoh(\P^N_k))$, then $Li^*G'$ generates $D(\Qcoh(\P^N_k))$.
\end{cor}

\proof Since $G'$ generates, we can get $\oplus_{j=0}^N \mathcal{O}_{\P^N_k}(j)$ and hence $G$. \qed

Let us recall another notion of generation.

\begin{defn}
 Given a idempotent-closed triangulated category $\mathcal{T}$. We say that a sub-category
$\mathcal{S}$ classically generates $\mathcal{T}$ if the smallest idempotent-closed triangulated subcategory containing $\mathcal{S}$ is $\mathcal{T}$ itself.
\end{defn}

\begin{cor}
 If $X$ is smooth, any bounded complex of coherent sheaves that classically strongly generates $D^b(\Coh(\P^N_k))$ pulls back to a classical generator of $D^b(\Coh(X))$.
\end{cor}

\proof This follows immediately from the previous result and theorem $2.1$ of \cite{Nee92}, but we shall give a direct proof. Since $X$ is smooth, every coherent sheaf admits a finite locally-free resolution so it is enough to show that we can get all locally-free coherent sheaves. Let $F$ be locally-free and consider $G$ as above. Then, we get each $\mathcal{O}_X(j)$ from summands of $G$ and a finite number of triangles. By assumption, we can get $G$ from our generator and hence we have the twisting powers. Now resolve $F$
\begin{gather*}
 \cdots \ra \mathcal{O}_X(-r_l)^{n_l} \ra \cdots \ra \mathcal{O}_X(-r_1)^{n_1} \ra F \ra 0
\end{gather*}
Let $d$ be the dimension of $X$. Truncate the above sequences at the $d+1$-st step and let $D$ denote the cokernel. It is also locally-free and of finite rank. From Grothendieck's theorem, $\Ext^n(F,D) = H^n(F^{\vee} \otimes D) = 0$ for $n > d$. Thus, in the triangle,
\begin{center}
\leavevmode
\begin{xy}
 (-10,-10)*+{D[d+1]}="a"; (0,5)*+{0\ra \mathcal{O}_X(-r_d)^{n_d} \ra \cdots \ra \mathcal{O}_X(-r_1)^{n_1} \ra 0}="b"; (10,-10)*+{F}="c"; {\ar@{->} "a";"b"}; {\ar@{->} "b";"c"}; {\ar@{->}^{\ \ \ \ \ \ a} "c";"a"}
\end{xy}
\end{center}
where $a$ is zero. The result then follows from the following lemma.

\begin{lem}
 Given a triangle
\begin{center}
\leavevmode
\begin{xy}
 (-10,10)*+{X}="a"; (10,10)*+{Y}="b"; (0,-5)*+{Z}="c"; {\ar@{->}^0 "a";"b"}; {\ar@{->} "b";"c"}; {\ar@{->} "c";"a"}
\end{xy}
\end{center}
 then $Z \cong Y \oplus X[1]$.
\end{lem}

\proof Given any $X'$ we have a short exact sequence
\begin{gather*}
 0 \ra \Hom(X',Y) \ra \Hom(X',Z) \ra \Hom(X',X[1]) \ra 0
\end{gather*}
Plugging in $X' = Y$ we see that there is a map $Z \ra Y$ that is left inverse to $Y \ra Z$. This splits the short exact sequence and gives the isomorphism \qed

\bigskip

\bibliographystyle{nhamsplain}
\bibliography{refup}

\end{document}